\documentclass[draft]{amsart}

\usepackage{graphics}
\usepackage[all]{xy}

\usepackage{amssymb,amsmath}
\usepackage{psfrag}

\usepackage{amsfonts}
\usepackage{amscd}

\usepackage{graphicx}
\usepackage{epstopdf}
\usepackage{tikz}
\usetikzlibrary{arrows,automata}

\usepackage{graphicx,epsfig,float}

\usepackage{url}

\makeatletter
\def\url@leostyle{%
  \@ifundefined{selectfont}{\def\UrlFont{\sf}}{\def\UrlFont{\small\ttfamily}}}
\makeatother
\newtheorem{thm}{Theorem}[section]
\newtheorem{lem}[thm]{Lemma}
\newtheorem{prop}[thm]{Proposition}
\newtheorem{cor}[thm]{Corollary}

\newtheorem{clm}[thm]{Claim}

\newtheorem{defn}[thm]{Definition}

\newtheorem{nrmk}[thm]{Remark}
\newtheorem{expl}[thm]{Example}

\newcommand{\pf}{{\bf Proof. }}
\newcommand{\ra}{\rangle}

\newcommand{\la}{\langle}

\def\Star{\mathrm{Star}}
\def\st{\mathrm{st}}

\renewcommand{\bar}{\overline}

\newcommand{\R}{\mbox{${\mathcal R}$}}

\newcommand{\id}{\mathrm{id}}

\newcommand{\bb}[1]{\ensuremath{\mathbb{#1}}}

\newcommand{\cal}[1]{\ensuremath{\mathcal{#1}}}

\begin{document}

\title {Coverings by open cells}

\author {M\'{a}rio J. Edmundo}

\address{ Universidade Aberta and CMAF Universidade de Lisboa\\
Av. Prof. Gama Pinto 2\\
1649-003 Lisboa, Portugal}

\email{edmundo@cii.fc.ul.pt}

\author {Pantelis  E. Eleftheriou}

\address{ Department of Pure Mathematics, University of Waterloo\\
200 University Ave West\\
N2L 3G1 Waterloo, Ontario, Canada}

\email{pelefthe@uwaterloo.ca}

\author{Luca Prelli}

\address{ CMAF Universidade de Lisboa\\
Av. Prof. Gama Pinto 2\\
1649-003 Lisboa, Portugal}

\email{lprelli@math.unipd.it}

\date{\today}
\thanks{The first author was supported by Funda\c{c}\~ao para a Ci\^encia e a Tecnologia, Financiamento Base 2008 - ISFL/1/209. The second author was supported by the Funda\c{c}\~ao para a Ci\^encia e a Tecnologia grant SFRH/BPD/35000/2007. The third author was supported by Marie Curie grant PIEF-GA-2010-272021. This work is part of the FCT project PTDC/MAT/101740/2008.\newline
 {\it Keywords and phrases:} O-minimal structures, open cells, semi-bounded structures.}

\subjclass[2010]{03C64}

\begin{abstract}
We prove that in a semi-bounded o-minimal expansion of an ordered group every non-empty open definable set is a finite union of open cells.
\end{abstract}

\maketitle

\begin{section}{Introduction}\label{section introduction}
We fix an arbitrary o-minimal expansion  ${\mathcal R}=\la R,<, +, 0, \dots\ra$ of an ordered group. Recall that by \cite{e1} $\R$ is semi-bounded if it has no poles; that is, in $\R$ there is no definable bijection between a bounded and an unbounded interval. See \cite{e1} for other characterizations of semi-boundeness.  In this note we prove  the following theorem.

\begin{thm}\label{thm}
If $\R$ is semi-bounded, then every non-empty open definable set is a finite union of  open cells.
\end{thm}

As explained in \cite[Subsection 2.1]{pet-sbd2}, there are three possibilities for  an arbitrary o-minimal expansion  ${\mathcal R}=\la R,<, +, 0, \dots\ra$ of an ordered group:

\begin{enumerate}
\item[(A)] $\R$ is linear (that is, its first-order theory ${\rm Th}(\R)$ is linear (\cite{lp})). In this case by \cite{lp}, there exists ${\mathcal S} \equiv \R$ with ${\mathcal S}$ a reduct of an ordered vector space $\mathcal{V}=\langle V, <, +, 0, \{d\}_{d\in D}\rangle$ over an ordered division ring $D$ (with the same addition and linear  ordering the underlying group of ${\mathcal S}$).

\item[(B)] $\R$ is not linear. In this case, the theory of every interval in $\R$ with the induced structure is not linear and so no interval in $\R$ is elementarily equivalent to a reduct of an interval in an ordered vector space (\cite{lp})). Therefore, by the Trichotomy theorem (\cite[Theorem 1.2]{pest-tri}), a real closed field whose ordering agrees with that of $\R$ is definable on some interval $(-e, e).$ There are now two sub-cases to consider:
\begin{itemize}
\item[(B1)] $\R$ is semi-bounded.
\item [(B2)] $\R$ is not semi-bounded. In this case, one can endow the whole structure $\R$ with a definable real closed field. Indeed, let $\sigma :(a, b)\to (c,+\infty )$ be a pole in $\R$; that is, a definable bijection (with say, $\lim _{t\to b}\sigma (t)=+\infty$). Without loss of generality, and using translations, we may assume that $a=c=0$ and $b<e.$ But then, being inside a real closed field, the intervals $(0,e)$ and $(0,b)$  are in definable bijection and  so $(0,e)$ and  $(0, +\infty )$ are in definable bijection. Now it is easy to get a real closed field on the whole of $\R.$
\end{itemize}

\end{enumerate}


A version of Theorem \ref{thm} in the field  case (B2) was proved by  Wilkie in \cite{w}, for bounded open definable subsets. There are simple examples that show that in this case the boundedness assumption is required. On the other hand, a version of Theorem \ref{thm} in the linear case (A) was proved by Andrews in \cite{and}. Here we generalize these two results to  the semi-bounded non-linear case. Moreover, we also prove a stronger result in the linear case, which we state next. For the notion of `linear decomposition' and `star', see Section \ref{sec-linear} below. For the notion of  `stratification', see \cite[Chapter 4, (1.11)]{vdd}.  By Lemma \ref{goodstrat}, Corollary \ref{contain}  and Proposition \ref{prop}  below, we have:

\begin{thm}\label{thm-str}
Assume that $\R=\la R, <,  0, +,  \{\lambda\}_{\lambda\in D}\ra$ is an ordered vector space over an ordered division ring $D.$ Let $\mathcal{D}$ be a linear decomposition  of $R^n$. Then there is  decomposition $\mathcal{C}$ of $R^n$ that refines $\mathcal{D}$, such that for every $C\in \mathcal{C}$, the star of $C$ is an open (usual) cell.
Moreover, $\mathcal C$ is a stratification of $R^n$.
\end{thm}


An important example of a semi-bounded, non-linear o-minimal structure is the expansion $\cal B$ of the real ordered vector space $\bb R_{vect}=\la \bb{R}, <, +, 0, \{d\}_{d\in \bb{R}}\ra$  by all bounded semi-algebraic sets. Every bounded interval in $\cal B$  admits the structure of a definable real closed field. For example, the field structure on $(-1, 1)$ induced from $\bb R$ via the semi-algebraic bijection $x\mapsto \frac{x}{\sqrt{1+x^2}}$ is definable in \cal B. By \cite{pss, mpp, pet-sbd1}, $\cal B$ is the unique structure that lies strictly between $\bb R_{vect}$ and the real field. The situation becomes significantly more subtle when $\mathcal R$ is non-archimedean, and the study of definable sets and groups in the general semi-bounded setting has recently regained a lot of interest (\cite{ee-sbg, el-sbd, ep-sbd1, ep-sbd2, pet-sbd2}).

We expect that our main theorem on coverings by open cells (Theorem \ref{thm}) will find numerous applications in the theory of locally definable manifolds in o-minimal structures. Some of those are exhibited in \cite{eep}. As stated in that reference, a strengthened result of coverings would yield further applications. We state the desired result here as a Conjecture:\vskip.2cm

\noindent\textbf{Conjecture.} Every definable set is a finite union of relatively open definable subsets which are definably simply connected.\vskip.2cm

\emph{Structure of the paper.} Section \ref{sec-linear} contains the stratification result (Theorem \ref{thm-str}) for the linear case. Section \ref{sec-sbd} contains the covering by open cells (Theorem \ref{thm}) for the semi-bounded non-linear case.
\end{section}

\medskip

\emph{Notation.} We recall the standard   notation  for graphs and ``generalized cylinders'' of definable maps.

\begin{itemize}
\item
If $f:X\to R$ is a definable map, we denote by $\Gamma(f)$ the graph of $f$.

\item
If $f, g:X \to R$ are definable maps or the constant maps $-\infty $ and $+\infty$ on $X$ with $f(x)<g(x)$ for all $x\in X$, we write $f<g$ and set:
\begin{itemize}
\item[]
$(f,g)_X=\{(x, y)\in X\times R: f(x)<y<g(x)\}$;

\item[]
$[f,g)_X=\{(x, y)\in X\times R: f(x)\leq y<g(x)\}$;

\item[]
$(f,g]_X=\{(x, y)\in X\times R: f(x)< y\leq g(x)\}$;

\item[]
$[f,g]_X=\{(x, y)\in X\times R: f(x)\leq y\leq g(x)\}$.
\end{itemize}
We also use the same notation for functions $f, g:Y\to R$ whose domain $Y$ contains $X$ and whose restrictions on $X$ are as above.
\end{itemize}

\emph{Acknowledgements.} We wish to thank the referee for many helpful comments.

\begin{section}{The linear case}\label{sec-linear}



We assume in this section that $\R=\la R, <,  0, +,  \{\lambda\}_{\lambda\in D}\ra$ is an ordered vector space over an ordered division ring $D$. For basic properties on such o-minimal structures we refer the reader to \cite[Chapter 1, Section 7]{vdd}.\\

A   function $f: R^n\to  R$  of the form $f(x_1, \ldots, x_n)=\lambda_1 x_1+\ldots+\lambda_n x_n+a$, where  $\lambda_i\in D$ and $a\in R$, is called  \emph{linear} (or  \emph{affine)}.  For a definable set $X\subseteq  R^n$, we denote by $L(X)$ the set of restrictions on $X$ of linear functions and by $L_\infty(X)$ the set $L(X)\cup\{\pm\infty\},$ where we regard $-\infty$ and $+\infty$ as constant functions on $X$. The functions from $L(X)$ are called linear functions on $X$. Clearly, if two linear functions have the same restrictions on $X$ then their restrictions on $cl (X)$ are equal as well. \\

We define {\it linear cells in $R^n$} inductively as follows:
\begin{itemize}
 \item
 a linear cell in $R$ is either a singleton subset of $R$, or an open interval with endpoints in $R\cup\{\pm\infty\}$,
 \item
 a linear cell in $R^{n+1}$ is a set of the form $\Gamma(f)$, for some $f\in L(X)$, or $(f,g)_X$, for some $f,g\in L_\infty(X)$, $f<g$, where $X$ is a linear cell in $R^n$.
\end{itemize}
In either case, $X$ is called \emph{the domain} of the defined cell. \\

We refer the reader to \cite[Chapter 3, (2.10)]{vdd} for the definition of a \emph{decomposition} of $R^n$. A \emph{linear decomposition of $R^n$} is then a decomposition $\mathcal{C}$ of $R^n$ such that each $B\in \mathcal{C}$ is a linear cell. The following can be proved similarly to \cite[Chapter 3, (2.11)]{vdd}.

\begin{thm}[Linear CDT]\label{thm LCDT}
$\,\,$
\begin{enumerate}
\item
Given any definable sets $A_1, \dots, A_k\subseteq  R^n$, there is a linear decomposition $\mathcal{C}$ of $R^n$ that partitions each $A_i$.

\item
Given a definable function $f:A\to  R$, there is a linear decomposition $\mathcal{C}$ of $R^n$ that partitions $A$ such that the restriction $f_{|B}$ to each $B\in\mathcal{C}$ with $B\subseteq  A$ is linear.
\end{enumerate}
\end{thm}


\begin{defn}\label{star}
{\em
Let $\mathcal{C}$ be a linear decomposition of $R^n$ and $X$ a subset of $R^n$.  Denote
$$\Star _{\mathcal{C}}(X)=\{D\in \mathcal{C}: X \cap cl(D)\ne\emptyset\}.$$
The \emph{star of $X$ with respect to $\mathcal{C}$}, denoted by $\st _\mathcal{C}(X)$, is then
$$\st _\mathcal{C}(X)=\bigcup \Star_\mathcal{C}(X).$$
We just write $\Star(X)$ and $\st(X)$ if  $\mathcal{C}$ is clear from the context.
}
\end{defn}


In what follows, if $k>0$, then $\pi:R^{k+1}\to   R^k$ denotes the usual projection map onto the first $k$-coordinates, and if $\mathcal{C}$ is a linear decomposition of $R^{k+1}$, then $\pi(\mathcal{C})$ denotes the linear decomposition
$\{\pi(C):C\in\mathcal{C}\}$ of $R^k$.

\begin{lem}\label{st}
Let $\mathcal{C}$ be a linear decomposition of $R^n$ and $X$ a subset of $R^n$. Then:

(i) If $n>1$, then $\Star_{\pi(\mathcal{C})}(\pi(X))=\pi(\Star_{\mathcal{C}}(X))$.

(ii) If $X$ is  an open  union of cells in $\mathcal{C}$, and $C\in \mathcal{C}$ with $C\subseteq  X$, then $\st(C)\subseteq  X$.
\end{lem}

\pf
(i) $\subseteq $. Let $D\in \Star(\pi(X))$.  Since $\pi$ is open, for any open set $U$ containing $X$, $\pi(U)$ is an open set containing $\pi(X)$. Thus $D\cap \pi(U)\ne \emptyset$, which implies $\pi^{-1}(D)\cap U\ne\emptyset$. Hence, by the definition of linear decomposition, there is some $D^\prime\in \Star(X)$ such that $\pi(D^\prime)=D$.

$\supseteq$. Let $D\in \Star(X)$. For any open set $U$ containing $\pi(X)$, $\pi^{-1}(U)$ is an open neighborhood of $X$. Therefore $\pi^{-1}(U)\cap D \ne\emptyset$, and $U\cap \pi(D)\ne\emptyset$. Hence $\pi(D)$ belongs to $\Star(\pi(X))$.

(ii) Since $X$ is open, for every $B\in \Star(C)$, $B\cap X\ne\emptyset$, and hence $B\subseteq  X$.
\qed\vskip.2cm

One would expect that $\st_\mathcal{C}(X)$ is an open set. However, the following example shows that this is not the case.

\begin{expl}\label{expl star}
{\em
Consider points $a_{-1}<a_0<a_1<a_2<a_3$ in $R$ and let $\mathcal{C}$ be a linear decomposition of $R^2$ that contains the following cells: $(a_{-1}, a_0)\times (a_0,a_2),\, (a_0, a_1)\times (a_0,a_2),\, \{a_0\}\times (a_0, a_1),\, \{a_0\}\times (a_1, a_3)$ and the point $(a_0,a_1)$. Then the star of the point $(a_0,a_1)$ is the union of the above cells, which is not open.
}
\end{expl}

Below we define a special kind of a linear decomposition $\mathcal{C}$ of $R^n$ that remedies the above problem. In fact, such a $\mathcal{C}$  will give us that every $\st_\mathcal{C}(X)$ with $X\in {\mathcal C}$ is an open (usual) cell (see Proposition \ref{prop} below). From this we obtain the version of Theorem \ref{thm} for the linear case (see Corollary \ref{cor cover open cells linear} below).\\



\begin{defn}\label{special}
{\em
A \emph{special linear decomposition of $R^n$} is a linear decomposition of $R^n$ defined by induction on $n$ as follows. Any linear decomposition of $R$ is special. A
linear decomposition $\mathcal{C}$ of $R^{k+1}$, $k>0$,  is special if:
\begin{itemize}
\item
$\pi (\mathcal{C})$ is a special linear decomposition of $R^k$;

\item
for any cells $\Gamma (h_{|A})$ and $(f,g)_B$ in $\mathcal{C}$, there is no $c\in cl(A)\cap cl(B)$ such that $f(c)<h(c)<g(c).$
\end{itemize}

}
\end{defn}

Before providing the nice consequences of special linear decompositions, we prove that they always exist.

\begin{lem}\label{goodstrat}
For any linear decomposition $\mathcal{D}$ of $R^n$, there is a special linear decomposition $\mathcal{C}$ of $R^n$ that refines $\mathcal{D}$ (that is, every linear cell in $\mathcal{D}$ is a union of linear cells in $\mathcal{C}$).
\end{lem}

\pf
By induction on $n$. For $n=1$, take $\mathcal{C}=\mathcal{D}$. Now assume that $n=k+1$ and the lemma holds for $k>0$. Let $\mathcal{D}$ be a linear decomposition of $R^{k+1}$. Choose a finite collection $\mathcal{F}$ of linear maps $f:R^k\to R$ such that any linear map that appears in the definition of any linear cell from $\mathcal{D}$ is a restriction of a map from  $\mathcal{F}$.
Now set
\[
\mathcal{G}=\{\Gamma(f)\cap\Gamma(g): f, g\in \mathcal{F}\} \text{ and } \mathcal{G}^\prime = \{\pi (A): A\in \mathcal{G}\} \cup \pi ({\mathcal D}).
\]
Clearly, $\mathcal{G}^\prime$ is a finite collection of definable subsets of $R^k$. By the linear CDT and the inductive hypothesis, there is a special
linear decomposition $\mathcal{C}^\prime$ of $R^k$ that partitions each member of $\mathcal{G}^\prime$.\vskip.2cm

\begin{clm}\label{clm goodstrat}
For any $f,g\in \mathcal{F}$, either $f< g$ or $f= g$ or $f> g$ on any $B\in \mathcal{C}'$.
\end{clm}

Let $B\in \mathcal{C}'$ and let $A=\Gamma(f)\cap \Gamma(g)$. Since $\pi(A)$ is a union of members of $\mathcal{C}'$, we have either $B\subseteq \pi (A)$ or $B\cap \pi (A)=\emptyset .$ In the first case $f=g$ on $B.$ In the second case, $B$ is a disjoint union of the open definable subsets $\{b\in B: f(b)<g(b)\}$ and $\{b\in B: g(b)<f(b)\}.$ Since $B$ is definably connected, one of the two sets is equal to $B$.
\qed
\vskip.2cm

Let $\mathcal{C}$ be the linear cell decomposition of $R^{k+1}$ with $\pi (\mathcal{C})=\mathcal{C}'$ such that for any $B\in \mathcal{C}'$ the set of cells in $\mathcal{C}$ with domain $B$ is defined by all functions from $\mathcal{F}.$ Since $\mathcal{C}'$ refines $\pi (\mathcal{D})$, the choice of $\mathcal{F}$  and Claim \ref{clm goodstrat} imply that $\mathcal{C}$ refines $\mathcal{D}.$

To conclude we need to show that $\mathcal{C}$ is special. Let $(f,g)_B\in \mathcal{C}.$ Then $f,g \in \mathcal{F}$ and for any $h\in \mathcal{F}$ we have on $B$ either $h< f$, or $h=f$, or $h= g$ or $h> g,$ and so either $h(c)\leq f(c)$ or $g(c)\leq h(c)$, for any $c\in cl(B).$ In particular,  for any $\Gamma (h_{|A})\in \mathcal{C}$ there is no $c\in cl(A)\cap cl(B)$ such that $f(c)<h(c)<g(c).$

\qed \\

We now aim towards Proposition \ref{prop} below. But before we will require several preliminary lemmas. \\

\begin{lem}\label{lem cl of cells linear}
Let $(R,<)$ be a dense linear order, $X\subseteq R^n$, and $\bar{X}=cl(X).$ Let  $f,g:\bar{X} \to R$ be continuous functions, and $C=(f,g)_X.$ Then
\begin{enumerate}
\item
$cl(\Gamma (f_{|X})=\Gamma (f);$
\item
$cl(C)=[f,g]_{\bar{X}}.$
\end{enumerate}
 In particular, $\pi (cl(C))=cl(\pi (C)).$

\end{lem}

\pf
(1) This is a special case of a general simple fact about continuous maps in Hausdorff topological spaces.

(2) Clearly, $[f,g]_{\bar{X}}=cl(C')$, where $C'=(f,g)_{\bar{X}}.$ So we need to show that $cl(C')=cl(C).$ Since $C\subseteq C'$, it is enough to show that $C'\subseteq cl(C)$. Let $(x,y)\in C'$ and let $U\times (a,b)$ be an open neighborhood of $(x,y)$ with $U$ an open neighborhood of $x$ and $a<y<b$. We may assume that $f(x)<a<y<b<g(x).$ Since $f$ and $g$ are continuous at $x$, there is an open $V$ with $x\in V\subseteq U$ such that $f(v)<a$ and $b<g(v)$ for all $v\in V.$ Since $x\in \bar{X}$, there is $v\in V\cap X$; and so, $(v,y)\in (U\times (a,b))\cap C.$ Thus $(x,y)\in cl(C)$ as required.

\qed \\

For $X\subseteq R^n$ a subset, $x=(x_1, \ldots , x_n)\in X$  and $\epsilon \in R^{>0}$ below we let $V_X(x, \epsilon )=\{u=(u_1,\ldots , u_n)\in X: |x_i-u_i|<\epsilon \,\,\textrm{for all}\,\, i\}$ (the {\it $\epsilon$-neighborhood of $x$ in $X$}).\\

\begin{lem}\label{lem linear x and eps}
If $X\subseteq R^n$ is a linear cell and $x\in X$ then there is $\epsilon \in R^{>0}$ such that $2x-y\in X$ for all $y\in V_X(x, \epsilon ).$
\end{lem}

\pf
We prove the result by induction on $n.$ Let $n=1.$ If $X$ is a singleton, take any $\epsilon \in R^{>0}.$ If $X$ is an open interval, take any $\epsilon \in R^{>0}$ such that $V_X(x, \epsilon) \subseteq X.$

Suppose that the  result holds for $n$ and we prove it for $n+1.$

Let $X=\Gamma (f_{|Z})$, and $x=(z, f(z))$ where $z\in Z$. By induction, there is $\epsilon \in R^{>0}$ such that $2z-u\in Z$ for all $u\in V_Z(z,\epsilon).$  By linearity  of $f$, it follows that, if  $y=(u, f(u)) \in V_X(x, \epsilon )$, then   $2x-y=(2z-u, 2f(z)-f(u))=(2z-u, f(2z-u))\in X$ since $u\in V_Z(z,\epsilon).$

Let $X=(f,g)_Z$, and $x=(z, v)$, where $z\in Z$ and $f(z)<v<g(z).$ Fix $\delta '\in R^{>0}$ such that $(v-\delta ', v+\delta ')\subseteq (f(z), g(z))$ and by continuity of $f$ and $g$ fix $\delta \in R^{>0}$ such that $(v-\delta ', v+\delta ')\subseteq (f(u), g(u))$ for all $u\in V_Z(z, \delta ).$ By induction, there is $\epsilon \in R^{>0}$ such that $2z-u\in Z$ for all $u\in V_Z(z,\epsilon).$  Choose $\epsilon <\delta ', \delta .$ Let $y=(u, w)\in V_X(x, \epsilon)$. Then $u\in V_Z(z,\epsilon)$ and $2z-u\in Z$;  $w\in (v-\epsilon, v+\epsilon )\subseteq (f(u), g(u))$ and so $2v-w\in (v-\epsilon, v+\epsilon)\subseteq (f(u), g(u)).$  Therefore,  $2x-y=(2z-u,2v-w)\in X.$
\qed \\

The following consequence of  Lemma \ref{lem linear x and eps} will be useful below:

\begin{lem}\label{compare clm}
Let $\mathcal{C}$ be a special linear decomposition of $R^n$, $n>1$ and $D\in \mathcal{C}$ a linear cell of the form
\[
D=\Gamma(f),
\]
where $f\in L(B).$  If $A\subseteq cl(B)$ is any linear cell, then there is a linear cell $F\in {\mathcal C}$ of the form $F=(h,k)_B$ such that $f_{|A}=h_{|A}<k_{|A}.$
\end{lem}

\pf
Since $A\subseteq cl(B)$, for any $h$ and $k$ with $(h,k)_B\in {\mathcal C}$ we have $h_{|A}\leq k_{|A}.$ Therefore, there is a linear cell $(h,k)_B\in {\mathcal C}$ which is above $D=\Gamma (f)$ and is  such that $f_{|A}=h_{|A}\neq k_{|A}$ and $h_{|A}\leq k_{|A}.$ We show that $h_{|A}<k_{|A}.$

If $k=+\infty $ the claim holds. Assume that $k\neq +\infty $ and let $p=k-h.$ Then $p_{|A}\geq 0.$ We have to show that $p_{|A}>0$. If not let $a\in A$ be such that $p_{|A}(a)=0.$ By Lemma \ref{lem linear x and eps} here is $\epsilon \in R^{>0}$ such that $2a-b\in A$ for all $b\in V_A(a, \epsilon ).$ Since $\{b\in A:p_{|A}(b)>0\}$ is an open definable subset of $A$ which is non empty (because $h_{|A}\neq k_{|A}$) and $A$ is definably connected, $\{b\in A: p_{|A}(b)=0\}$ is a closed non open definable subset of $A.$ So there is a $c\in V_A(a, \epsilon )$ such that $p_{|A}(c)>0.$ But then since $2a-c\in A$, we have $p_{|A}(2a-c)=2p_{|A}(a)-p_{|A}(c)=-p_{|A}(c)<0$ contradicting the fact $p_{|A}\geq 0.$
\qed\\

Below we also need the following remark.

\begin{nrmk}\label{nrmk convex}
{\em
Let $A\subseteq R^n$ be a subset. We say that $A$ is {\it convex} if for all $x,y\in A$ and for all $q\in {\mathbb Q}\cap [0,1]$ we have $qx+(1-q)y\in A.$ See \cite[Definition 3.1]{el-st}.

The following hold:
\begin{itemize}
\item[$\bullet$]
The intersection of two convex sets is convex.
\item[$\bullet $]
Every linear cell is convex.\\
\end{itemize}
}
\end{nrmk}

We are now ready to prove the main lemma for what follows below.

\begin{lem}\label{compare}
Let $\mathcal{C}$ be a special linear decomposition of $R^n$, $n>1$, $D, E\in \mathcal{C}$ two linear cells of the form
\[
D=\Gamma(f) \text{ and }\, E=\Gamma(g),
\]
where $f\in L(B)$, $g\in L(A)$, and $A\subseteq cl(B)$. Then:
\[
\text{$f_{| A}< g$ or $f_{| A}= g$ or $f_{| A}> g$}.
\]
\end{lem}

\pf
Assume not. Then there is $c\in A$ such that $f(c)=g(c);$ otherwise $A$ would the be disjoint union of the open definable subsets $\{x\in A: f(x)<g(x)\}$ and $\{x\in A: g(x)<f(x)\}$ contradicting the fact that $A$ is definably connected. Since $f_{|A}\neq g$, then there exists $d\in A$ such that $f(d)\ne g(d).$ We may assume that $f(d)<g(d).$ By Lemma \ref{compare clm}, there is a linear cell $F\in \mathcal{C}$ of the form $F=(h, k)_B$ such that $f_{| A}=h_{| A}< k_{| A}$. We next show that there is a point $e\in A$, such that $h(e)<g(e)<k(e)$ which contradicts the fact that ${\mathcal C}$ is special.

$\,$
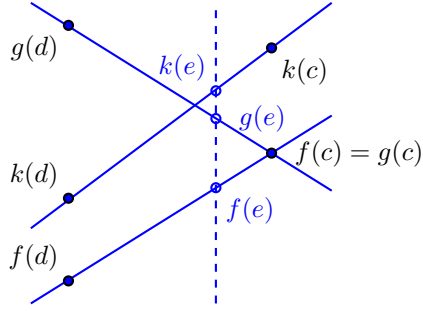
\begin{center}

\begin{tikzpicture}[domain=0:2]

\draw[thick, color= blue] (-2, -1) -- (2, 2);
\draw[thick, color= blue] (-2, -2) -- (2,0.5);
\draw[thick, color= blue] (-2, 2) -- (2,-0.5);

\draw[thick, color= blue, dashed] (0.46, -2) -- (0.46,2);

\draw[thick, fill= blue] (-1.5,1.7) circle (.6mm) node[below left] {$g(d)$};
\draw[thick, fill= blue]  (-1.5,-0.6) circle (.6mm) node[above left] {$k(d)$};
\draw[thick, fill= blue]  (-1.5,-1.7) circle (.6mm) node[above left] {$f(d)$};

\draw[thick, color= blue]  (0.46,0.83) circle (.6mm) node[above left] {$k(e)$};;
\draw[thick, color= blue]  (0.46,0.46) circle (.6mm) node[right ] {$\,\,\, g(e)$};;
\draw[thick, color= blue]  (0.46,-0.46) circle (.6mm) node[below right] {$f(e)$};;

\draw[thick, fill= blue] (1.2,1.4) circle (.6mm) node[below right] {$k(c)$};
\draw[thick, fill= blue] (1.2,0) circle (.6mm) node[right] {$\,\,\, f(c)=g(c)$};

\end{tikzpicture}

\end{center}

If $g(d)<k(d)$, then let $e=d$. So assume $k(d)\le g(d)$. We will choose $e$ to be ``between'' $c$ and $d$. We first see that there is $q_0\in (0,1]\cap \mathbb{Q}$, such that
\[
q_0 g(d)+(1-q_0)g(c)<q_0 k(d)+(1-q_0)k(c)
\]
Indeed, if not, then $k(c)\le g(c)$. But $g(c)=f(c)<k(c)$, a contradiction. On the other hand, since $f(d)<g(d)$ and $f(c)=g(c)$, we have that
for every $q\in (0,1]\cap \mathbb{Q}$,
\[
qf(d)+(1-q)f(c)<q g(d)+(1-q)g(c).
\]
Hence, if we let $e=q_0 d+(1-q_0) c$, then $e\in A$ (by Remark \ref{nrmk convex}) and we have $f(e)=h(e)<g(e)<k(e)$, proving our claim.


\qed \\

\begin{lem}\label{picl}
Let $\mathcal{C}$ be a special linear decomposition of $R^n$, $n>1$, and $D, E\in \mathcal{C}$ such that $D\cap cl(E)\ne\emptyset$. Then:
\[
\pi(D)\subseteq  cl(\pi(E))\, \Rightarrow \, D\subseteq  cl(E).
\]
\end{lem}

\pf
Let $A=\pi (D)$ and $B=\pi (E)$; so $A\subseteq cl(B).$ We have the following possibilities for $E$: (1) $ E=\Gamma (f_{|B})$ or (2) $E=(f,g)_B;$ and the following possibilities for $D$: (a) $D=\Gamma (h_{|A})$ or (b) $D=(h,k)_A.$

By Lemma \ref{lem cl of cells linear}, if (1) then $cl(E)=\Gamma (f_{|cl(B)})$, and if (2) then $cl(E)=[f,g]_{cl(B)}.$

Suppose (1). If (a), since $D\cap cl(E)\ne\emptyset$,  there is $a\in A$ with $f(a)=h(a)$ and so by Lemma \ref{compare}, $f_{|A}=h_{|A}$ and therefore $D\subseteq cl(E).$ On the other hand, case (b) under (1) cannot happen: as  ${\mathcal C}$ is special, there is no $a\in A$ such that $h(a)<f(a)<k(a)$ and so $D\cap cl(E)=\emptyset $ contradicting the assumption of the lemma.

Suppose (2).  If (a), since $D\cap cl(E)\ne\emptyset$,  there is $a\in A$ with $f(a)\leq h(a)\leq g(a).$ Since ${\mathcal C}$ is special, $f(a)=h(a)$ or $h(a)=g(a)$ and so by Lemma \ref{compare}, $f_{|A}=h_{|A}$ or $h_{|A}=g_{|A}$ and therefore $D\subseteq cl(E).$ If (b), there is $a\in A$ such that $[f(a), g(a)]$ intersects $(h(a), k(a)).$ Since ${\mathcal C}$ is special, we have $f(a)=h(a)$ and $g(a)=k(a)$ and so by Lemma \ref{compare}, $f_{|A}=h_{|A}$ or $g_{|A}=k_{|A}$ and therefore $D\subseteq cl(E).$



\qed \\

\begin{cor}\label{contain}
Let $\mathcal{C}$ be a special linear decomposition of $R^n$, $n>0$, and $D, E\in \mathcal{C}$ such that $D\cap cl(E)\ne\emptyset$. Then $D\subseteq  cl(E)$.

In particular, $\mathcal C$ is a stratification of $R^n$.
\end{cor}

\pf
The statement trivially holds if $D=E$, hence assume $D\ne E$. We work by induction on $n$. For $n=1$, the assumption $D\cap cl(E)\ne\emptyset$
implies that $E$ is an open interval and $D$ is one of its endpoints. So now assume $n>1$. Clearly, $\pi(D)\cap cl(\pi(E))\ne\emptyset$ (using Lemma \ref{lem cl of cells linear}), and
hence by inductive hypothesis, $\pi(D)\subseteq  cl(\pi(E))$. By Lemma \ref{picl}, $D\subseteq  cl(E)$.
\qed

\begin{lem}\label{starpi}
Let $\mathcal{C}$ be a special linear decomposition of $R^n$, $n>0$. Then, for any subset $X\subseteq  R^n$, $\st(X)$ is open.
\end{lem}

\pf
It suffices to show that $\st (X)\cap cl(E)=\emptyset $ for any $E\in \mathcal{C}$ with $\st(X)\cap E=\emptyset. $ Suppose this is not the case. Then  some $D\in \Star (X)$ meets $cl(E)$. Then by Corollary  \ref{contain}, $cl(E)$ contains $D$ and so $cl(D)$. As $X$ meets $cl(D)$, it meets $cl(E)$, and hence $E\subseteq \st(X)$, which is a contradiction.



\qed

\begin{prop}\label{prop}
Let $\mathcal{C}$ be a special linear decomposition of $R^n$, $n>0$, and $C\in \mathcal{C}$. Then $U=\st(C)$ is an open (usual) cell.
\end{prop}

\pf
By Lemma \ref{starpi}, $U$ is open. So it remains to prove that $U$ is a cell. Before that we need a few preliminaries.

Since $ \mathcal{C}$ is a linear decomposition, for every $B\in \Star(\pi(C))$, $\pi^{-1}(B)\cap U$ is
a union of linear cells in $\mathcal{C}$ which are either graphs of linear maps, or cylinders between linear maps, with domain $B$. By Lemma \ref{st}(i), $U\subseteq  \bigcup\{\pi^{-1}(B): B\in \Star(\pi(C))\}$, and hence
\[
U= \bigcup\{\pi^{-1}(B)\cap U: B\in \Star(\pi(C))\}.
\]
We claim that for every $B\in \Star(\pi(C))$,
\[
 \pi^{-1}(B)\cap U=(f_B, g_B)_B,
\]
for some $f_B, g_B\in L_\infty(B)$ with $f_B< g_B$.

Fix $B\in \Star(\pi(C)).$ Let $f_B$ be the bottom function with domain $B$ defining the bottom cell of $\pi^{-1} (B) \cap U$  and let $g_B$ be the top function with domain $B$ defining  the top cell of $\pi^{-1} (B) \cap U.$ (Recall that this latter set is a union of linear cells in $\mathcal{C}$ which are either graphs of linear maps, or cylinders between linear maps, with domain $B$).

\begin{clm}\label{clm cyl main}
If $C=(l,k)_P$ then $E=(f_B, g_B)_B$ is the unique cell in $\Star(C)$ such that $\pi (E)=B$ and $f_{B|P}=l$ and $g_{B|P}=k.$ In particular, $\pi^{-1}(B)\cap U=(f_B, g_B)_B.$
\end{clm}

 Suppose that $E$ is not a  cell in $\Star(C)$. Then there  are cells $\Gamma (h^1_{|B}), \ldots , \Gamma (h^m_{|B})$ in $\Star(C)$ such that $f_B<h^1_{|B}<\dots <h^m_{|B}<g_B.$ Since $B\in \Star(\pi(C))$ we have $P=\pi (C)\subseteq cl(B),$ and since $\Gamma (h^i_{|B})\in \Star (C)$ we have $C\subseteq cl(\Gamma (h^i_{|B})=\Gamma (h^i_{|cl(B)})$ (by Lemma \ref{lem cl of cells linear}). Hence $cl(C)= [l,k]_{cl(P)}\subseteq \Gamma (h^i_{|cl(B)})$ (by Lemma \ref{lem cl of cells linear}) and therefore, $l=(h^i _{|cl(B)})_{|P}=k$ which is absurd.

 By the choice of $f_B$ and $g_B$, $E=(f_B, g_B)_B$ is then the unique cell in $\Star(C)$ such that $\pi (E)=B.$ Since $C\subseteq cl (E)=[f_B,g_B]_{cl(B)}$  we have $cl(C)= [l,k]_{cl(P)}\subseteq  [f_B,g_B]_{cl(B)}.$ Since  ${\mathcal C}$ is special we must have $f_{B|P}=l$ and $g_{B|P}=k.$
 \qed\\

\begin{clm}\label{clm grp main}
If $C=\Gamma (l_{|P})$ then there  are cells $\Gamma (h^1_{|B}), \ldots , \Gamma (h^m_{|B})$ in $\Star(C)$ such that $f_B<h^1_{|B}<\dots <h^m_{|B}<g_B$ and $h^i _{B|P}=l.$ Moreover, $(h^i_{|B}, h^{i+1}_{|B})_B, (f_B, h^1_{|B})_B$ and $(h^m _{|B}, g_B)$ are all cells in $\Star(C).$ In particular, $\pi^{-1}(B)\cap U=(f_B, g_B)_B.$
\end{clm}

Let $<h^1_{|B}<\dots <h^m_{|B}$ be all the linear functions that appear in the definition of a linear cell of $\pi^{-1} (B) \cap U$. (Recall that $U=\st(C)$ and $\pi^{-1} (B) \cap U$ is a union of linear cells in $\mathcal{C}$ which are either graphs of linear maps, or cylinders between linear maps, with domain $B$).  Then $\Gamma (h^1_{|B}), \ldots , \Gamma (h^m_{|B})$ are cells in $\Star(C)$ and  $(h^i_{|B}, h^{i+1}_{|B})_B, (f_B, h^1_{|B})_B$ and $(h^m _{|B}, g_B)$ are all cells in $\Star(C).$

Since $P\subseteq B$ and $\Gamma (h^i_{|B})\in \Star (C)$, by Lemma \ref{compare} we must have $h^i _{B|P}=l.$
\qed \\

 We conclude the proof of the proposition by induction on $n$. If $n=1$, then $C$ is a point and $U$ is an open interval or $C$ is an open interval and $U=C$. Now assume that $n=k+1$ and the result holds for $k>0$.

Let $D=\st(\pi(C))$, $f=\bigcup_{B\in \Star(\pi(C))} f_B$ and $g=\bigcup_{B\in \Star(\pi(C))} g_B$. Then
\[
U=(f, g)_D.
\]
By inductive hypothesis, $D$ is a usual cell.
To show that $f, g$ are continuous, we need to show that for every $A,B\in \Star(\pi(C))$,
and $A\subseteq cl(B)$,
\[
f_{B|A}=f_A \text{ and } g_{B|A}=g_A.
\]
Indeed, for any $B, B'\in \Star (\pi (C))$, if $cl(B)\cap cl(B')\neq \emptyset ,$ then the intersection of $cl(B)\cap cl(B')$ with the domain of $f$ (resp. $g$) is a union of cells  $A\in \mathcal{C}$ such that  $A\subseteq cl(B)\cap cl(B')$ (by Corollary \ref{contain}) and $A\in \Star(\pi (C)).$


By Lemma \ref{compare}, there are 3 possibilities: (i) $f_{B|A}>f_A$, (ii) $f_{B|A}=f_A$, (iii) $f_{B|A}<f_A$.

If we assume (i) we get a contradiction since in that case $U$ is not open. Let us assume (iii).

If $C=\Gamma (l_{|P})$, then since ${\mathcal C}$ is special, by Lemma \ref{compare},  and using the notation of Claim \ref{clm grp main}, we have $(h^1_{|B})_{|A}\leq f_A.$ Since $P\subseteq cl(A)\subseteq cl(B),$ we have $l=(h^1_A)_{|P}>(f_A)_{|P}\geq (h^1_{|B})_{|P}=l$, which is absurd.

If $C=(l,k)_P$, then by Claim \ref{clm cyl main}, $F=(f_A,g_A)_A$ and $E=(f_B, g_B)$ are in $\Star(C)$ and $l=(f_A)_{|P}=(f_B)_{|P}$ and $k=(g_A)_{|P}=(g_B)_{|P}.$ Since ${\mathcal C}$ is special, by Lemma \ref{compare},  we have $(g_{B})_{|A}\leq f_A.$ Since $P\subseteq cl(A),$ we have $(g_B)_{|P}\leq (f_A)_{|P}$. Hence, if $(x,y)\in C$, then $l(x)<y<k(x)=(g_B)_{|P}(x)\leq (f_A)_{|P}(x)$ and so $(x,y)\notin [f_A, g_A]_{cl(A)}=cl(F)$. So $C\not \subseteq cl(F)$ which is absurd.
\qed\\

\begin{cor}\label{cor cover open cells linear}
If $\R=( R, <,  0, +,  \{\lambda\}_{\lambda\in D})$ is an ordered vector space over an ordered division ring $D$, then every non-empty open definable set is a finite union of open cells.

\end{cor}

\pf
Let $X\subseteq R^n$ be an open definable subset and take $\mathcal{C}$ a special linear decomposition of $R^n$ that partitions $X$. By Lemma \ref{st}(ii),
\[
X=\bigcup_{C\in\mathcal{C},\, C\subseteq  X} \st(C).
\]
Then apply Proposition \ref{prop}.
\qed
\end{section}

\begin{section}{The semi-bounded non-linear case}\label{sec-sbd}

We assume in this section that  $\R$ is semi-bounded and non-linear.  So, as we saw in the Introduction,  there exists a definable real closed field $\la I,0_I, 1_I,+_I, \cdot _I,<_I\ra$ on some interval $I\subseteq R$ which, without loss of generality,  can be assumed to be of the form  $I=(-e,e)$, $0_I=0$ and $<_I $ is the restriction of $<$ to $I$. 
Here we will use the existence of this  ``short'' definable real closed field to adapt Wilkie's proof  (\cite{w})  in o-minimal expansions of real closed fields. \\




In the next  lemmas the semi-boundedness assumption of $\R$ is not required.

\begin{lem}[\cite{w}, Lemma 1]\label{lem w lem1}
 Let $C$ be a cell in $R^n$. Then there exists an open cell $D$ in $R^n$ with $C\subseteq D$ and a definable retraction $H:D\to  C$ (that is, a continuous map such that $H_{|C}=\id _C$).
\end{lem}


\begin{lem}\label{lem nl cos}
Let $C$ be a cell in $R^n$. Suppose that $h:C\to  R$ is a continuous definable map and let $U$ be an open definable subset of $R^{n+1}$. Suppose further that  $\Gamma (h)\subseteq U$. Then there exist definable maps $f,g:C\to  R$ and cells $C_1,\dots ,C_m\subseteq C$ such that:
\begin{enumerate}
\item
$f <h <g$;
\item
$C=C_1\cup \cdots \cup C_m$;
\item
for each $i$, $f_{|C_i}$ and $g_{|C_i}$ are continuous;
\item
for each $i$, $\Gamma (h_{|C_i})\subseteq [f_{|C_i},g_{|C_i}]_{C_i}\subseteq U$.
\end{enumerate}
\end{lem}

\pf
Since $U$ is open and  $\Gamma (h)\subseteq U$, by definable choice (\cite[Chapter 6, (1.2)]{vdd}  there exists definable maps $f,g:C\to  R$ such that $f <h <g$ and $[f ,g]_{C}\subseteq U.$ By cell decomposition, there are cells  $C_1,\dots ,C_m\subseteq C$  covering $C$ such that for each $i$, $f_{|C_i}$ and $g_{|C_i}$ are continuous. Now the rest is clear.
\qed \\

The following is also needed:

\begin{lem}\label{lem nl cos2}
Let $C$ be a cell in $R^n$. Suppose that $f,g:C\to  R$  are continuous definable maps such that $f<g$ and let $V, W\subseteq U$ be  open definable subsets of $R^{n+1}$. Suppose further that  $(f,g)_C\subseteq U$, $\Gamma (f)\subseteq V$ and $\Gamma (g)\subseteq W$. Then there exist definable maps $f',g':C\to  R$ and cells $C_1,\dots ,C_m\subseteq C$ such that:
\begin{enumerate}
\item
$C=C_1\cup \cdots \cup C_m$;
\item
for each $i$, $f'_{|C_i}$ and $g'_{|C_i}$ are continuous;
\item
$f <f'<g' <g$;
\item
for each $i$, $\Gamma (f'_{|C_i})\subseteq V$ and $\Gamma (g'_{|C_i})\subseteq W$;
\item
for each $i$, $(f'_{|C_i},g_{|C_i})_{C_i}\subseteq U$, $(f _{|C_i},g'_{|C_i})_{C_i}\subseteq U$ and $[f'_{|C_i},g'_{|C_i}]_{C_i}\subseteq U$.
\end{enumerate}
\end{lem}

\pf
Since $(f,g)_C\subseteq U$, $\Gamma (f)\subseteq V$ and $\Gamma (g)\subseteq W$ and  $V, W\subseteq U$ be  open definable subsets of $R^{n+1}$, by definable choice (\cite[Chapter 6, (1.2)]{vdd}  there exists definable maps $f',g':C\to  R$ such that
\begin{enumerate}
\item
$f <f'<g' <g$;
\item
$\Gamma (f')\subseteq V$ and $\Gamma (g')\subseteq W$;
\item
$(f',g)_{C}\subseteq U$, $(f ,g')_{C}\subseteq U$ and $[f', g']_{C}\subseteq U$.
\end{enumerate}
By cell decomposition, there are cells  $C_1,\dots ,C_m\subseteq C$  covering $C$ such that for each $i$, $f_{|C_i}$ and $g_{|C_i}$ are continuous. Now the rest is clear.
\qed \\

Below we let
$$d^{(n)}(\bar{x},\bar{y})=\max \{|x_i-y_i| 1\leq i\leq n\}$$
denote the standard  distance in $R^n$ (where we denote by $\bar{z}=(z_1, \ldots , z_n)$ the elements of $R^n$). This distance is a continuous definable function (by \cite[Chapter 6 (1.4)]{vdd}). Moreover, if $B\subseteq R^n$ is a nonempty definable subset and $\bar{a}\in R^n$, then
$$d^n(\bar{a}, B)=\inf \{d^n(\bar{a},\bar{x}):\bar{x}\in B\}$$
is well defined (by (\cite[Chapter 1 (3.3)]{vdd})) and $d^n(\bar{a}, B)=0$ if and only if $\bar{a}\in cl(B)$ (the if part of this equivalence is immediate and for the only if part one can use the curve selection (\cite[Chapter 6 (1.5)]{vdd})).\\

Let  $\pi :R^{n+1}\to  R^n$ be the projection onto the first $n$ coordinates.  We say that an open definable subset $U$ of $R^{n+1}$  {\it has  $I$-short height} if for every $\bar{x}\in \pi (U)$ we have $${\rm sup}\{|t-s|:t,s\in U_{\bar{x}}\}\in I$$
where $U_{\bar{x}}=\{y\in R: (\bar{x},y)\in U\}.$\\

We now prove the  analogue of  \cite[Lemma 2]{w} for open definable subsets with $I$-short height. The argument of the proof is similar, one just has to observe that the field operations are used in Wilkie's proof in a uniform way and only along fibers. Since in our case our fibers are $I$-short, such field operations, in the  field $I$, can also be used in exactly the same way.

For completeness we include the details of the proof but at the end  we follow a more constructive argument  suggested to us  by Oleg Belagradek. For that we need the following observations which are true in arbitrary o-minimal expansions of ordered groups:

\begin{nrmk}\label{nrmk oleg1}
{\em
If $\theta :[a, b]\to [c,d]$ is a continuous, definable, strictly decreasing function, $\theta (a)=d$ and $\theta (b)=c$, then $\theta $ is bijective.

Indeed, as $\theta $ is definable and continuous, $\theta ([a,b])$ is definable, closed, and bounded by \cite[Chapter 6 (1.10)]{vdd}, and hence it a finite union of closed intervals and singletons, by o-minimality. Since $\theta $ is strictly decreasing, $\theta ([a,b])$ is densely ordered, and so is a closed interval, which must be $[c,d].$
}
\end{nrmk}

\begin{nrmk}\label{nrmk oleg2}
{\em
Let $V\subseteq R^n$ be an open definable subset, and let $\{\theta _{\bar{x}}:[a_{\bar{x}}, b_{\bar{x}}]\to [c_{\bar{x}}, d_{\bar{x}}] \}_{\bar{x}\in V}$ be a uniformly definable family of strictly decreasing functions with $\theta _{\bar{x}}(a_{\bar{x}})=d_{\bar{x}}$ and $\theta _{\bar{x}}(b_{\bar{x}})=c_{\bar{x}}$. (So by the previous remark all $\theta _{\bar{x}}$'s are bijective). Suppose that all $a_{\bar{x}}, b_{\bar{x}}, c_{\bar{x}}, d_{\bar{x}}$ are continuous functions in $\bar{x}$, and moreover, the map
$$\{(\bar{x},y):\bar{x}\in V \,\,\textrm{and}\,\, y\in [a_{\bar{x}}, b_{\bar{x}}]\} \to R: (\bar{x},y)\mapsto \theta _{\bar{x}}(y)$$ is  continuous. Then the map
$$\gamma : \{(\bar{x},z):\bar{x}\in V \,\,\textrm{and}\,\, z\in [c_{\bar{x}}, d_{\bar{x}}]\} \to R: (\bar{x},y)\mapsto \theta _{\bar{x}}^{-1}(z)$$
is continuous.

Indeed, for each $\bar{x}$, let $\bar{\theta }_{\bar{x}} :R\to R$ be given by
\begin{equation*}
\bar{\theta }_{\bar{x}}(y)=
\begin{cases}
a_{\bar{x}} +d_{\bar{x}} -y \qquad \textrm{for} \,\,\, y<a_{\bar{x}}
\\
\theta _{\bar{x}}(y) \qquad \,\,\,\,\, \qquad \textrm{for} \,\,\, a_{\bar{x}}\leq y \leq b_{\bar{x}}
\\
 b_{\bar{x}}+ c_{\bar{x}} -y \,\, \qquad \textrm{for} \,\,\, y>b_{\bar{x}}
\end{cases}
\end{equation*}
Then $\{\bar{\theta } _{\bar{x}}: R\to R\}_{\bar{x}\in V}$ is a uniformly definable family of strictly decreasing functions such that $\bar{\theta }_{\bar{x}}$ extends $\theta _{\bar{x}}$ for all $\bar{x}\in V$, and $V\times R\to R: (\bar{x},y)\mapsto \bar{\theta }_{\bar{x}} (y)$  is a continuous function.

Now $\bar{\gamma } : V\times R\to R: (\bar{x},y)\mapsto \bar{\theta }_{\bar{x}}^{-1} (z)$  is also a continuous function since, for any $(a,b)\subseteq R,$
\begin{eqnarray*}
\bar{\gamma }^{-1}((a,b)) & = & \{(\bar{x}, z)\in V\times R: a<\bar{\theta }_{\bar{x}}^{-1} (z)<b\} \\
& = & \{(\bar{x}, z)\in V\times R: \bar{\theta }_{\bar{x}}(a)<z<\bar{\theta }_{\bar{x}}(b)\},
\end{eqnarray*}
which is open. Therefore, since $\bar{\gamma }$ extends $\gamma $, we have that $\gamma $ is also continuous as required.
}
\end{nrmk}

\begin{lem}\label{lem sb wilkie lem2 short}
Let $C$ be a cell in $R^n$. Suppose that $f,g:C\to  R$ are continuous definable maps such that $f<g$ and let $U$ be an open definable subset of $R^{n+1}$ with $I$-short height. Suppose further that $[f,g)_C\subseteq U$  (respectively $(f,g]_C\subseteq U$). Then there exists an open definable subset $V$ of $R^n$ and continuous definable maps $F,G:V\to  R$ such that:
\begin{enumerate}
\item
$C\subseteq V$;
\item
$F_{|C}=f$ and $\Gamma (F)\subseteq U$ (respectively $\Gamma (G)\subseteq U$);
\item
$G_{|C}=g$;
\item
$F<G$;
\item
for all $\bar{x}\in V$ and all $y\in R$ with $F(\bar{x})\leq y<G(\bar{x})$, (respectively $F(\bar{x})<y\leq G(\bar{x})$), $(\bar{x},y)\in U$.
\end{enumerate}
\end{lem}

 \pf
 We prove the unparenthesized statement, the parenthetical one being similar.

 Applying Lemma \ref{lem w lem1} we obtain an open cell $D$ in $R^n$, with $C\subseteq D$, and a continuous definable retraction $H:D\to  C$.

Let
$$V=\{\bar{x}\in D: d^{(n)}(\bar{x},H(\bar{x}))<d^{(n+1)}((\bar{x},f\circ H(\bar{x})), U^c)\},$$
where $U^c=R^{n+1}\setminus U.$ Clearly $V$ is open in $R^n$ and (1) holds since $\Gamma (f)\subseteq U.$ Putting $F=f\circ H_{|V}$  we see that (2) holds. Also note that  for all $\bar{x}\in V$,  $F(\bar{x})<g\circ H(\bar{x})$ and
$$J_{\bar{x}}:=[0, g\circ H(\bar{x})-F(\bar{x}))\subseteq \{t\in R_{\geq 0}:F(\bar{x})+t\in U_{\bar{x}}\}\subseteq  I$$
since $U$ has $I$-short height.

By o-minimality and the fact that $\Gamma (F)\subseteq U$, there are well defined definable maps $z_0:V\to  I$ and $y_0:V\to  R$ given by

$$z_0(\bar{x})={\rm sup}\{t\in J_{\bar{x}}: [F(\bar{x}), F(\bar{x})+t)\subseteq U_{\bar{x}}\} $$
and $$y_0(\bar{x})=F(\bar{x})+z_0(\bar{x}).$$
Now observe that $y_0:V\to  R$ satisfies the conditions (3), (4) and (5) for $G$ ((3) is satisfied because $(f,g)_C\subseteq U$, by hypothesis, and $f=F_{|C}$), but maybe $y_0$ is not continuous. Thus we need to find a continuous definable map $G:V\to  R$ such that $F<G\leq y_0$ and $G_{|C}=y_0.$

Consider the definable set
$$S=\{(\bar{x},y)\in R^{n+1}:\bar{x}\in V\,\,{\rm and}\,\,F(\bar{x})\leq y\leq g\circ H(\bar{x})\}$$
and the definable continuous maps $\theta _1, \theta _2:S\to  I$ given by
$$\theta _1(\bar{x},y)=1_I - _I(y-F(\bar{x}))\cdot _I(g\circ H(\bar{x})-F(\bar{x}))^{- _I1_I}$$
where $1_I$ is the neutral element for the multiplication $\cdot _I$, $- _I$ is the diference and $ \,\, ^{- _I1_I}$ is inversion in the field $I$, and,
$$\theta _2(\bar{x},y)={\rm inf}\{d^{(n+1)}((\bar{x},t),U^c):F(\bar{x})\leq t\leq y\}.$$
Note that since $U$ has $I$-short height we do have $\theta _1(S)\subseteq I$ and   $\theta _2(S)\subseteq I.$

Fix $\bar{x}\in V$. Then the continuous definable map $(\theta _1 \cdot _I \theta _2)(\bar{x},-)$ decreases monotonically and strictly from $d^{(n+1)}((\bar{x},F(\bar{x})),U^c)$ to $0_I=0$ on $[F(\bar{x}),y_0(\bar{x})]$ and is identically $0_I=0$ on $[y_0(\bar{x}),g\circ H(\bar{x})].$

For $\bar{x}\in V$ let
\[
a_{\bar{x}}=F(\bar{x}), \,\, b_{\bar{x}}=y_0(\bar{x}), \,\, c_{\bar{x}}=0, \,\, d_{\bar{x}}=d^{(n+1)}((\bar{x},F(\bar{x})),U^c),
\]
and $\theta _{\bar{x}}(-)=(\theta _1 \cdot _I \theta _2)(\bar{x},-)_{|}: [a_{\bar{x}}, b_{\bar{x}}]\to  [c_{\bar{x}}, d_{\bar{x}}].$ Then by Remark \ref{nrmk oleg2},
$$G:V\to R:\bar{x}\mapsto \theta _{\bar{x}}^{-1}(d^{(n)}(\bar{x},H(\bar{x})))$$
is a continuous definable function. Moreover,  $a_{\bar{x}}< G(\bar{x})\leq b_{\bar{x}}$ for all $\bar{x}\in V.$ In fact, if not then $a_{\bar{x}}=G(\bar{x})$ and we obtain $(\theta _1\cdot _I \theta _2)(\bar{x},G(\bar{x}))=d^{(n+1)}((\bar{x},F(\bar{x})), U^c)$ contradicting the fact that $d^{(n)}(\bar{x},H(\bar{x}))<d^{(n+1)}((\bar{x},F(\bar{x})), U^c).$ We also have $G(\bar{x})=b_{\bar{x}}$ for all $\bar{x}\in C.$ Therefore, $G$ satisfies (3), (4) and (5) as required.




 \qed \\

We need one more  lemma:

\begin{lem}\label{lem sb wilkie lem2 closed}
Let $C$ be a cell in $R^n$. Suppose that $f,g:C\to  R$ are continuous definable maps such that  $f<g$  and let $U$ be an open definable subset of $R^{n+1}$. Suppose further that $[f,g]_C\subseteq U$. Then there exists an open definable subset $W$ of $R^n$ and continuous definable maps $F,G:W\to  R$ such that:
\begin{enumerate}
\item
$C\subseteq W$;
\item
$F_{|C}=f$ and $\Gamma (F)\subseteq U$;
\item
$G_{|C}=g$ and $\Gamma (G)\subseteq U$;
\item
$F<G$;
\item
for all $\bar{x}\in W$ and all $y\in R$ with $F(\bar{x})\leq y\leq G(\bar{x})$, $(\bar{x},y)\in U$.
\end{enumerate}
\end{lem}

\pf
Applying Lemma \ref{lem w lem1} we obtain an open cell $D$ in $R^n$, with $C\subseteq D$, and a continuous definable retraction $H:D\to  C$.

Let  $W'$ be the intersection of
$$\{\bar{x}\in D: d^{(n)}(\bar{x},H(\bar{x}))<d^{(n+1)}((\bar{x},f\circ H(\bar{x})), U^c)\}$$
and
$$\{\bar{x}\in D: d^{(n)}(\bar{x},H(\bar{x}))<d^{(n+1)}((\bar{x},g\circ H(\bar{x})), U^c)\}$$
where $U^c=R^{n+1}\setminus U.$ Clearly $W'$ is open in $R^n$ and (1) holds for $W'$ since $\Gamma (f), \Gamma (g) \subseteq U.$ Also (2) and (3)  hold for $f\circ H_{|W'}$ and $g\circ H_{|W'}$. Also note that  for all $\bar{x}\in W'$,  $f\circ H_{|W'}(\bar{x})<g\circ H_{|W'}(\bar{x})$ so (4) holds for $f\circ H_{|W'}$ and $g\circ H_{|W'}$.

Let $B=[f\circ H_{|W'}, g\circ H_{|W'}]_{|W'}\setminus U$ where
{\small
$$[f\circ H_{|W'}, g\circ H_{|W'}]_{|W'}=\{(\bar{x},y)\in W'\times R: y\in [f\circ H_{|W'}(\bar{x}), g\circ H_{|W'}(\bar{x})]\},$$
}
and let $$W=W'\setminus \bar{\pi (B)}.$$  Clearly $W$ is open. We now show that $C\subseteq W$, verifying in this way (1). Suppose not and let  $c=(c_1, \ldots , c_n)\in C$ be such that $c\in \bar{\pi (B)}$. Let $\epsilon >0$ be such that $E=\Pi _{i=1}^n[c_i-\epsilon, c_i+\epsilon]\subseteq W'$. By definable choice there is a definable map $\alpha :(0,\epsilon )\to  \pi (B)\cap E$ such that ${\rm lim}_{t\to 0^+}\alpha (t)=c.$ By replacing $\epsilon $ we may assume that $\alpha $ is continuous. Again by definable choice, we see that there exists a definable map $\beta :(0,\epsilon )\to  B\cap [f\circ H_{|E}, g\circ H_{|E}]_{|E}$ such that $\pi \circ \beta =\alpha$. By replacing $\epsilon $ we may assume that $\beta $ is continuous. Since the definable set $B\cap [f\circ H_{|E}, g\circ H_{|E}]_{|E}$  is closed and, by \cite[Proposition 3.1 (3)]{e1},  $\beta ((0,\epsilon ))$ is bounded,
the limit ${\rm lim}_{t\rightarrow  0^+}\beta (t)$ exists  in this set. If $d$ is this limit, then $\pi (d)=c$ since $\pi \circ \beta =\alpha.$ So $d\in [f\circ H_{|W'}(c), g\circ H_{|W'}(c)]\cap B\neq \emptyset $ contradicting the fact that $ [f\circ H_{|W'}(c), g\circ H_{|W'}(c)]=[f(c),g(c)]\subseteq U.$

If we put $F=f\circ H_{|W}$ and $G=g\circ H_{|W}$ we see that (2), (3) and (4) hold. On the other hand, if $\bar{x}\in W$ and $y\in R$ are such that $F(\bar{x})\leq y\leq G(\bar{x})$ and, by absurd, $(\bar{x},y)\notin U$, then $(\bar{x},y)\in B$ and so $\bar{x}\in \pi (B)\subseteq \bar{\pi (B)}$ contradicting the fact that $\bar{x}\notin \bar{\pi (B)}.$ Thus (5) also holds.
\qed \\

Combining Lemmas \ref{lem sb wilkie lem2 short} and \ref{lem sb wilkie lem2 closed} we obtain:

\begin{lem}\label{lem sb wilkie lem2 general}
Let $C$ be a cell in $R^n$. Suppose that $f,g:C\to  R$ are continuous definable maps such that $f<g$ and let $U$ be an open definable subset of $R^{n+1}$. Suppose further that $[f,g)_C\subseteq U$  (respectively $(f, g]_C\subseteq U$). Then there exists a cell decomposition $C_1,\dots ,C_l$ of $C$ and for each $i=1,\dots ,l$ there is an open definable subset $V_i$ of $R^n$ and continuous definable maps $F_i,G_i:V_i\to  R$ such that:
\begin{enumerate}
\item
$C_i\subseteq V_i$;
\item
$F_{i|C_i}=f_{|C_i}$ and $\Gamma  (F_i)\subseteq U$ (respectively $\Gamma (G_i)\subseteq U$);
\item
$G_{i|C_i}=g_{|C_i}$;
\item
$F_i<G_i$;
\item
for all $\bar{x}\in V_i$ and all $y\in R$ with $F_i(\bar{x})\leq y<G_i(\bar{x})$, (respectively $F_i(\bar{x})<y\leq G_i(\bar{x})$), $(\bar{x},y)\in U$.
\end{enumerate}
\end{lem}

 \pf
 We prove the unparenthesized statement, the parenthetical one being similar.

Let $H:D\to C$ be as in  Lemma \ref{lem w lem1}. Choose $\epsilon \in I$ such that $2\epsilon \in I$ and put

$$U_f= U\cap ( (f\circ H)-\epsilon, (f\circ H)+\epsilon )_{D} $$
and
$$U_g= U\cap ( (g\circ H)-\epsilon, (g\circ H)+\epsilon )_{D}.$$
Then clearly $U_f$ and $U_g$  are open definable  subsets  of $U$ with $I$-short height. For example,  if $(\bar{x}, y)\in U_f$, then $(f\circ H)(\bar{x})-\epsilon < y < (f\circ H)(\bar{x})+\epsilon.$

Since  $(f,g)_C\subseteq U$, $\Gamma (f)\subseteq U_f$ and $\Gamma (g)\subseteq U_g$, by Lemma \ref {lem nl cos2},  there exist definable maps $f',g':C\to  R$ and cells $C_1,\dots ,C_m\subseteq C$ such that:
\begin{enumerate}
\item
$C=C_1\cup \cdots \cup C_m$;
\item
for each $i$, $f'_{|C_i}$ and $g'_{|C_i}$ are continuous;
\item
$f <f'<g' <g$;
\item
for each $i$, $\Gamma (f'_{|C_i})\subseteq U_f$ and $\Gamma (g'_{|C_i})\subseteq U_g$;
\item
for each $i$, $(f'_{|C_i},g_{|C_i})_{C_i}\subseteq U$, $(f _{|C_i},g'_{|C_i})_{C_i}\subseteq U$ and $[f'_{|C_i},g'_{|C_i}]_{C_i}\subseteq U$.
\end{enumerate}

Fix $i=1,\dots ,m$. Then we can apply Lemma \ref{lem sb wilkie lem2 short} to the data $(U_f, f_{|C_i}, f'_{|C_i})$ and obtain the data $(V_f,F_1 ,F'_1)$ satisfying  (1) to (5) of that lemma. Similarly, we can apply Lemma \ref{lem sb wilkie lem2 short} to the data $(U_g, g'_{|C_i}, g_{|C_i})$ and obtain the data $(V_g,G'_1,G_1)$ satisfying  (1) to (5) of that lemma. On the other hand, we can apply Lemma \ref{lem sb wilkie lem2 closed} to the data $(U, f'_{|C_i}, g'_{|C_i})$ and obtain the data $(W,F',G')$ satisfying  (1) to (5) of that lemma .

Take $V_i=V_f\cap V_g\cap W$ and set $F=F_{1|V_i}$, $G=G_{1|C_i}.$ Then clearly (1) to (5) hold.
 \qed \\

The following is also required:

\begin{lem}\label{lem sb wilkie lem3 infty}
Let $C$ be a cell in $R^n$. Suppose that $k:C\to  R$ is a continuous definable map  and let $U$ be an open definable subset of $R^{n+1}$. Suppose further that $[k,+\infty )_C\subseteq U$  (respectively $(-\infty ,k ]_C\subseteq U$). Then there exists an open definable subset $W$ of $R^n$ and a continuous definable map $K:W\to  R$ such that:
\begin{enumerate}
\item
$C\subseteq W$;
\item
$K_{|C}=k$ and $\Gamma (K)\subseteq U$;
\item
for all $\bar{x}\in W$ and all $y\in R$ with $K(\bar{x})\leq y$ (respectively $y\leq K(\bar{x})$), $(\bar{x},y)\in U$.
\end{enumerate}
\end{lem}

\pf
We prove the unparenthesized statement, the parenthetical one being similar.

Applying Lemma \ref{lem w lem1} we obtain an open cell $D$ in $R^n$, with $C\subseteq D$, and a continuous definable retraction $H:D\to  C$.


Let
$$W'=\{\bar{x}\in D: d^{(n)}(\bar{x},H(\bar{x}))<d^{(n+1)}((\bar{x},k\circ H(\bar{x})), U^c)\}$$
where $U^c=R^{n+1}\setminus U.$ Clearly $W'$ is open in $R^n$ and (1) holds for $W'$ since $\Gamma (k)\subseteq U$. Also (2)   holds for $k\circ H_{|W'}$.

Let $B=[k\circ H_{|W'}, +\infty )_{|W'}\setminus U$ where
{\small
$$[k\circ H_{|W'}, +\infty )_{|W'}=\{(\bar{x},y)\in W'\times R: k\circ H_{|W'}(\bar{x})\leq y\},$$
}
and let $$W=W'\setminus \bar{\pi (B)}.$$  Clearly $W$ is open. We now show that $C\subseteq W$, verifying in this way (1). Suppose not and let  $c=(c_1, \ldots , c_n)\in C$ be such that $c\in \bar{\pi (B)}$. Let $\epsilon >0$ be such that $E=\Pi _{i=1}^n[c_i-\epsilon, c_i+\epsilon]\subseteq W'$.  By definable choice there is a definable map $\alpha :(0,\epsilon )\to  \pi (B)\cap E$ such that ${\rm lim}_{t\to 0^+}\alpha (t)=c.$ By replacing $\epsilon $ we may assume that $\alpha $ is continuous. Again by definable choice, we see that there exists a definable map $\beta :(0,\epsilon )\to  B\cap [k\circ H_{|E}, +\infty )_{|E}$ such that $\pi \circ \beta =\alpha$. By replacing $\epsilon $ we may assume that $\beta $ is continuous. Since the definable set $B\cap [k\circ H_{|E}, +\infty )_{|E}$  is closed and, by \cite[Proposition 3.1 (3)]{e1},  $\beta ((0,\epsilon ))$ is bounded,
the limit ${\rm lim}_{t\rightarrow  0^+}\beta (t)$ exists  in this set. If $d$ is this limit, then $\pi (d)=c$ since $\pi \circ \beta =\alpha.$ So $d\in [k\circ H_{|W'}(c), +\infty )\cap B\neq \emptyset $ contradicting the fact that $ [k\circ H_{|W'}(c), +\infty )=[k(c), +\infty )\subseteq U.$

If we put $K=k\circ H_{|W}$  we see that (2) holds. On the other hand, if $\bar{x}\in W$ and $y\in R$ are such that $K(\bar{x})\leq y$ and, by absurd, $(\bar{x},y)\notin U$, then $(\bar{x},y)\in B$ and so $\bar{x}\in \pi (B)\subseteq \bar{\pi (B)}$ contradicting the fact that $\bar{x}\notin \bar{\pi (B)}.$ Thus (3) also holds.
\qed

\begin{cor}\label{cor cover open cells sb non-linear}
If $\R$ is  a semi-bounded non-linear o-minimal expansion of an ordered group, then every non-empty open definable set is a finite union of  open cells.
\end{cor}

\pf
This is done by induction on the dimension of the open definable set. For dimension one this is clear. Let $U$ be an open definable subset of $R^{n+1}$. Let ${\mathcal D}$ be a  cell decomposition of $R^{n+1}$ partitioning  $U$. Clearly it is enough to show that each  cell $D\in {\mathcal D}$ with $D\subseteq U$ can  be covered by finitely many open  cells (in $R^{n+1}$) each of which is contained in $U$.

Case A: $D=(f,g)_C$ for some  cell $C$ in $R^n$ and continuous definable   maps $f,g:C\to  R$ such that  $f<g.$

Let $f'=\frac{2f+g}{3}$ and $g'=\frac{f+2g}{3}.$ Then $f',g':C\to  R$ are continuous definable maps such that
\begin{itemize}
\item
$f<f'<g'<g$;
 \item
 $\Gamma (f')\subseteq U$ and  $\Gamma (g')\subseteq U$;
\item
$(f', g)_C\subseteq U$ and $(f, g')_C\subseteq U$.
 \end{itemize}
Now apply Lemma \ref{lem sb wilkie lem2 general} to the data $(C,U, f, g')$ and obtain the data $(C_i,V_i,F_i ,G'_i)$ with $i=1,\dots ,l$ satisfying  (1) to (5) of that lemma.  By  the inductive hypothesis there exists a finite collection ${\mathcal A}_i$ of open cell in $R^n$ contained in $V_i$ which cover $V_i$.  By (4) and (5) of  Lemma \ref{lem sb wilkie lem2 general}, for each  $A\in {\mathcal A}_i,$ $(F_{i|A},G'_{i|A})_A$ is an open cell in $R^{n+1}$ contained in $U$, and by (1), (2) and (3) of that lemma, $(f_{|C_i},g'_{|C_i})_{C_i}\subseteq \cup \{(F_{i|A},G'_{i|A})_A:A\in {\mathcal A}_i\}.$ Thus  $(f,g')_{C}\subseteq \cup \{(F_{i|A},G'_{i|A})_A:A\in {\mathcal A}_i\,\,{\rm and}\,\,i=1,\dots ,l\}.$

Similarly, apply  Lemma \ref{lem sb wilkie lem2 general} to the data $(C,U, f', g)$ (the parenthetical statement there) to see that $(f',g)_C$ can be covered by finitely many open cells in $R^{n+1}$ each of which is contained in $U$. Hence the same is true for $(f,g)_C=(f,g')_C\cup (f',g)_C.$

Case B: $D=\Gamma (h)$ for some continuous definable  map $h:C\to  R$ where $C$ is a  cell in $R^n$. This case reduces to Case A above by Lemma  \ref{lem nl cos}.

Case C: $D=(k, +\infty )_C$ (respectively $D=(-\infty , k)_C$) for some  cell $C$ in $R^n$ and continuous definable   map $k:C\to  R$.

Then we can apply Lemma \ref{lem sb wilkie lem3 infty} to the data $(C,U, k)$ and obtain the data $(C,W, K)$ satisfying  (1) to (3) of that lemma.  By  the inductive hypothesis there exists a finite collection ${\mathcal A}$ of open cell in $R^n$ contained in $W$ which cover $W$.  By  (3) of  Lemma \ref{lem sb wilkie lem3 infty}, for each  $A\in {\mathcal A},$ $(K_{|A},+\infty )_A$ is an open cell in $R^{n+1}$ contained in $U$, and by (1) and (2) of that lemma, $(k_{|C},+\infty )_{C}\subseteq \cup \{(K_{|A},+\infty )_A:A\in {\mathcal A}\}.$

Similarly for the case $D=(-\infty , k)_C$.
\qed






\end{section}

\end{document}